\documentclass{article}               
\usepackage{color}
\usepackage{graphicx}

\usepackage{hyperref}
\hypersetup{ 
colorlinks=true, %colorise les liens 
breaklinks=true, %permet le retour à la ligne dans les liens trop longs 
urlcolor= blue, %couleur des hyperliens 
linkcolor= blue, %couleur des liens internes 
bookmarksopen=true, %si les signets Acrobat sont créés, 
bookmarksnumbered=true,
pdftitle={Hybrid control for low-regular nonlinear systems: application to an embedded control for an electric vehicle}, %informations apparaissant dans 
pdfauthor={Thomas Chambrion and Gilles Millerioux}, %dans les informations du document 
pdfsubject={This note presents an embedded automatic control strategy for a low 
 consumption vehicle equipped with an ``on/off'' engine with  low cost/ low consumption constraints for the control device.  
},
pdfkeywords={Hybrid control, Adaptative control, Optimal control, Robustness, Automotive.} %sous Acrobat. 
}

  \usepackage{amsmath}
\usepackage{amssymb}

\usepackage{url}

\newtheorem{thm}{Theorem}
\newtheorem{proposition}[thm]{Proposition}
\newtheorem{lemma}[thm]{Lemma}
\newtheorem{remark}{Remark}

\newtheorem{assumption}{Assumption}
\newtheorem{problem}{Problem}

\begin{document}

\title{Hybrid control for low-regular nonlinear systems: application to an embedded control for an electric vehicle
}

\author{Thomas Chambrion\footnote{
Universit\'e de Lorraine, IECL, UMR 7502,  54506 Vandoeuvre-l\`es-Nancy, France} 
\footnote{
CNRS, UMR 7502, 54506 Vandoeuvre-l\`es-Nancy, France} 
\footnote{Inria, 54600 Villers, France, 
\texttt{Thomas.Chambrion@univ-lorraine.fr}}
\and               
Gilles Millerioux\footnote{
Universit\'e de Lorraine, CRAN, UMR 7039, 
54506 Vandoeuvre-l\`es-Nancy, France,
\texttt{Gilles.Millerioux@univ-lorraine.fr}}}

\maketitle

\begin{abstract}                          
 This note presents an embedded automatic control strategy for a low 
 consumption vehicle equipped with an ``on/off'' engine. 
The main difficulties are the hybrid nature of the dynamics, the non 
smoothness of the dynamics of each mode, the uncertain environment, 
the fast changing dynamics,  and  low cost/ low consumption 
constraints for the control device.  
Human drivers of such vehicles frequently use an oscillating strategy, 
letting the velocity evolve between fixed lower and upper bounds.
We present a general justification of this very simple and efficient strategy, that happens to be 
optimal for autonomous dynamics,
robust and easily adaptable for real-time control strategy. Effective implementation in a competition  prototype involved 
in low-consumption races shows that automatic velocity control achieves
performances comparable with the results of trained human drivers. Major 
advantages of automatic control
 are improved robustness and safety. 
The total average power consumption for the control device is less than 
10 mW. 
\end{abstract}
\setcounter{tocdepth}{2}
\tableofcontents

\newpage 
\section{Introduction}

%\subsection{Context}\label{SEC_context}

\begin{figure}
\centering
\includegraphics[width=5in]{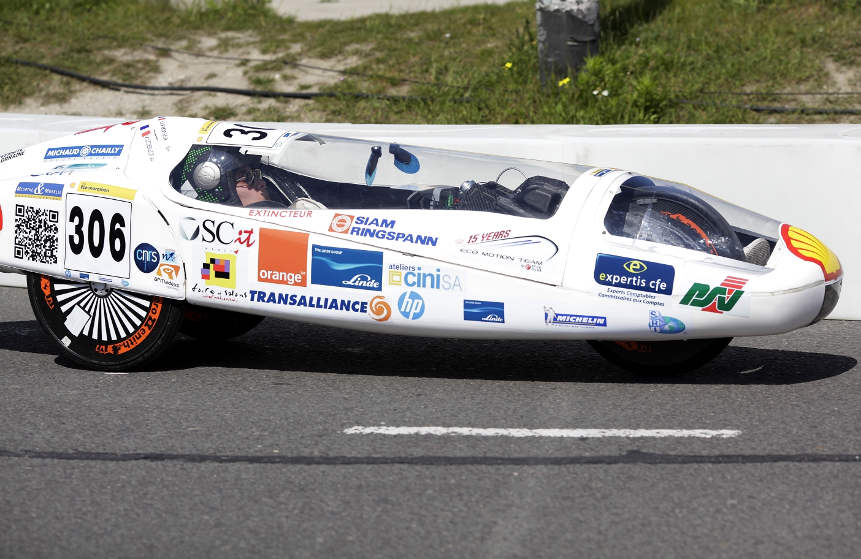}
\caption{Test run of Vir'Volt 2 in Rotterdam (May 2014).}
\label{fig_Essai}
\end{figure}

The
European Shell Eco-Marathon\textsuperscript{TM} brings together over 200 teams
from high schools and universities from all over Europe,
in a race involving ecological and economical vehicles. The
principle of the race is to go through a given distance
in a limited time and with the lower energy consumption.
%The Ecole Sup\'erieure des Sciences et Technologies de
%l'Ing\'enieur de Nancy (ESSTIN, France), has been involved in the
%European Shell Eco-marathon race since 15 years in the
%respective category fuel cell and battery. 
The aim of this note is to expose the design and the implementation of  an embedded automatic control of 
the speed on the prototype Vir'Volt 2 (see Figure~\ref{fig_Essai}), built by the students of \'Ecole Sup\'erieure des Sciences et Technologies de Nancy (ESSTIN) in France that 
participated to the edition 2014 of the European 
Shell Eco-Marathon.

A crucial feature of the prototype is the hybrid nature of the control. Indeed, the DC 
electric motor  is either switched on or switched off.
The energy consumption of the vehicle is reduced to a roughly constant residual 
 consumption of the electronics  (a few milliwatts) when the engine is off and goes up to 
 160 W  when the engine  is on. When the engine is off, the transmission 
 parts (gears, chains,...) are decoupled from the wheels and rapidly stops
 (the kinetic energy 
 is dissipated by friction and Joule effect in the engine).  
  When the engine is switched on, the transmission parts (gears, chains, ..) that were 
  at rest have to speed up in order to reach the wheels velocity, before the coupling allows   
  energy to  
  be transmitted from the engine to the wheels. The difference between the kinetic energy of the freely 
  evolving transmission parts and the transmissions parts rotating at engine speed is far from being negligible (about 
  10 J), and has to be paid each time the motor is switched on.

 %A few additional solar cells (0.17 m$^2$, peak power around 40 W) are allowed to 
%furnish energy all along the race, for   
%up to 20$\%$ of the total energy budget, without being incorporated in the official 
%consumption record.   

%\subsection{Former results}
Clearly, the issue enters the context of control in switched systems.
  Switching systems are widely involved in industrial processes \cite{871312}  and have been a major field of interest in automatic control  \cite{654885}\cite{liberz03} for years.
Optimal control of such switched systems is now pretty well understood for linear dynamics 
\cite{riedinger1999linear,Sun2005181,5208242,1039806}. For non linear dynamics, the case without switching cost has raised considerable attention, see for instance 
\cite{827981,1259455,MR2157228,1272479,6544639,1649936}. The optimization  problem is classically tackled with the Pontriaguin 
Maximum Principle (PMP), see for instance \cite{MR2062547} for an intrinsic formulation.
Non smooth extensions of the PMP have been developed  to include the case
of  non smooth trajectories or  nonzero switching costs \cite{832814,832811}. 

Two difficulties prevent the direct use of these results in our case.

The first difficulty arises when modeling the dynamics. In the standard frame of Lipschitz continuous dynamics, a vehicle with switched off engine on a flat road will never be at rest because of the principle of non-intersection of solutions. A more realistic modeling (where cars can actually stop) requires to consider less regular dynamics, which usually give rise to non-trivial technical subtleties. This explain that this approach is rare in practice (with few exceptions, such as \cite{898692} for instance).

 The second difficulty is the heavy computational burden. Even in the standard Lipschitz continuous framework, with especially adapted algorithm 
 \cite{5765434,5208312,1576861}, 
 some dozen of seconds  of computations are needed
 with a powerful desktop computer 
to compute the
 optimal strategy of a 39 minutes run, and any perturbation (traffic, weather change,...) requires a new 
computation. Following the rules of the competition, this computation has to be done or stored in the 
vehicle (not off board). Even for races where the electric energy needed to do the 
computation is not incorporated in the total energy consumption, the battery weight has bad 
influence on the total consumption. Online integration of the PMP requires 
a strong multi-core processor (about 50 W on average, about 100 kJ for the race). Table of optimal trajectories may be pre-computed offline 
 and stored inboard. The main drawback of this approach is the cost of 
 construction (pre-computation) of such a table - this includes a good knowledge of the geometry of the track and a reliable meteorological date base. This pre-computation process has to be redone for every different race. 
    
Since global optimization seems out of reach, one may consider finite horizon optimization (as in Model Predictive Control). Even in this case, the computational burden is not to be neglected: a small 
linux computer with an average power of $5$ W requires about $10$ kJ for the race, to be compared with the 
100 kJ needed to move the vehicle. Moreover, the choice of the horizon is usually difficult because it impacts on both feasibility \cite{Bemporad99,CamachoLIB04,LimonAUT08} and stability issues (see \cite{GilbertTanTAC91,KolmanovskyGilbert98,BlanchiniBook,manriqueala2016} for computational aspects with invariant set-based approaches).

%\textbf{REF,REF} and the stability is not guaranteed \textbf{REF,REF}.

%\subsection{Contribution}
Any proposed strategy to adapt the prototype speed will be an imperfect trade-off between
 robustness, computational tractability and efficiency. In energy efficient prototypes competitions, human drivers usually have an oscillating strategy, letting the velocity evolve between a lower and an upper limit. Our main result is a rigorous proof that, with physically reasonable assumptions,  the long term optimal strategy for a stationary problem (wind and track slope remain constant)  has exactly this pattern.   
Taking advantage of the particular structure of these solutions, 
this paper gives an explicit control strategy %(see Section \ref{SEC_strategy}) 
for the general problem (variable wind and slope) to adapt  the finite
horizon of the optimization scheme, with the following features:
\begin{itemize}
\item Very low computational burden, implementable on low-cost embedded micro-controller for real-time automatic control;
\item Very low power consumption  (less than 10 mW in average, 
about 15 J for the race);
\item Real time adaptability with respect to changes in dynamics, weather or traffic conditions;
\item Optimality of the control law  if the race conditions (wind and slope) are constant.
\end{itemize}
The result applies also in the case (important in practice) where the dynamics is not Lipschitz continuous.

% \subsection{Content of the paper}
The layout of the paper is the following.
In Section~\ref{SEC_Autonome_simplified}, the time independent model of the dynamics is presented and a finite-time horizon is considered. For this special case (autonomous 
dynamics corresponding to constant slope of the track, wind strength and direction),
 we prove (Proposition \ref{PRO_structure_minimale})  that the optimal driving strategies exhibit
  a very simple behavior, where the speed oscillates, most of the time, between a maximum and a minimum without intermediate local extrema. In Section \ref{SEC_simplified_periodic}, we consider a long term version of the problem studied in Section \ref{SEC_Autonome_simplified}, that is infinite time horizon. Once again, the optimal strategy exhibits a very simple periodic pattern. In Section \ref{SEC_Framework}, a non autonomous model for the dynamics is considered by taking into account possible external disturbances like variable slope, wind, \ldots. Some preliminary well-posedness mathematical results are introduced. Next, we use the strong geometric structure of the optimal strategies in the autonomous case (the finite horizon of the optimization scheme is the period of the optimal control computed in Section \ref{SEC_Autonome_simplified}). Stability  and robustness 
   analyses completes this section. 
Finally, Section \ref{SEC_implementation} is devoted to the results obtainef from a real-life implementation on the prototype during the competition.   

\section{Autonomous dynamics in finite time}\label{SEC_Autonome_simplified}

\subsection{Statement of the optimal control problem}\label{SEC_problem_statement_simplified}
Let us consider the unperturbed dynamics of the vehicle given by:
\begin{equation}\label{EQ_dynamics_simplified}
\dot{x}= f(x,u)
\end{equation}
where $x$ is the real-valued speed of the prototype,  
 $u$ is the control, a piecewise constant function taking value in $\{0,1\}$ representing the state of the engine (0 is off, 1 is on), and $f:\mathbf{R} \times \{0,1\} \to \mathbf{R}$ is the dynamics. The initial condition $x(0)=x_0\geq 0$ corresponds to the speed at start.  \\

Hereafter, the total duration of the race will be denoted with $T$ ($T>0$), the length of the track will be denoted with $L$ ($L>0$) and the average speed will be denoted with $V=L/T$ ($V>0$).

Let $u:[0,T]\to \{0,1\}$ be piecewise constant control. For a real constant $T_1>0$, a real-valued function 
$x:[0,T_1]\to \mathbf{R}$ is a  \emph{solution} of 
(\ref{EQ_dynamics_simplified}) if $x$ is  continuous, piecewise differentiable and 
satisfies (\ref{EQ_dynamics_simplified}) for  almost every time $t$ in $[0,T_1]$.

Let us recall that a real-valued function $g:[0,T_1]\to \mathbf{R}$ is absolutely continuous if $g$ is differentiable almost everywhere in $[0,T_1]$
with integrable derivative $g'$ that satisfies $g(t)=g(0)+\int_{0}^t g'(s)\mathrm{d}s$ for every $t$ in $[0,T_1]$.  With every admissible control $u$ and absolutely continuous solution $x$ of (\ref{EQ_dynamics_simplified}), we associate the energy cost
\begin{equation}\label{EQ_def_cout}
C(u)=\int_0^T h\big(x(t),u(t)\big)\mathrm{d}t+ \alpha\,N(u)\
\end{equation}
where $h:[0,+\infty)\times \{0,1\} \to [0,+\infty)$ is a $C^1$ function accounting for the instantaneous power consummed by the vehicle, $\alpha$ is a 
positive real constant (representing the cost to switch the motor on) and $N$ denotes the number of times that $u$ jumps from $0$ to 1.

 The problem to be solved is the following
 \begin{problem}\label{PROB_general_autonomous}
 Find a piecewise constant control $u:[0,T]\to \{0,1\}$ such that the associated solution of (\ref{EQ_dynamics_simplified}) with given initial condition $x_0$ satisfies $\int_0^T x(t)\mathrm{d}t=L$ and  minimizes the cost (\ref{EQ_def_cout}).
 \end{problem}

To ensure that the problem is well-posed, we need some regularity assumptions on $f$.   

\begin{assumption}\label{ASS_maths_simplified}
\begin{enumerate}
\item For any $u$ in $\{0,1\}$, the one-variable function $x\mapsto f(x,u)$ is continuous\label{ASS_maths_simplified_continuity}.
\item For every $x_0$ in $\mathbf{R}$, for every $u$ in $\{0,1\}$, the Cauchy problem \eqref{EQ_dynamics_simplified} with initial condition $x(0)=x_0$ admits a unique solution in positive time\label{ASS_maths_simplified_uniqueness}.
\end{enumerate}
\end{assumption}

The existence (but not the uniqueness) of continuously differentiable solutions to the Cauchy problem (\ref{EQ_dynamics_simplified}) with initial condition $x_0$ at time $t_0$ is a consequence of Assumption \ref{ASS_maths_simplified}.\ref{ASS_maths_simplified_continuity} and Peano theorem (see for instance 
\cite{MR0171038}, Theorem 2.1 page 10).
The local uniqueness in positive times of the solutions (Assumption \ref{ASS_maths_simplified}.\ref{ASS_maths_simplified_uniqueness}) is needed
for physical reasons and ensures the continuity of the solutions with respect to the initial conditions (\cite{MR0171038}, Theorem 2.1 page 94). Notice that Assumption \ref{ASS_maths_simplified} is not enough in general to ensure neither the uniqueness in past times nor the differentiability of the solutions with respect to the initial conditions (\cite{MR0171038}, Theorem 6.1 page 104).

\subsection{Preliminaries}\label{SEC_preliminaries_reparam}
Let  $u:[t_1,t_2]\subseteq[0,T]\to \{0,1\}$ be 
constant on the interval $[t_1,t_2]$. Then, the solution $x:[t_1,t_2]\to 
\mathbf{R}$ of (\ref{EQ_dynamics_simplified}) with initial value $x(t_1)\subseteq[0,L]$ is of class $C^1$. If moreover $t\mapsto f\big(x(t),u(t)\big)$ does not vanish on $[t_1,t_2]$, then $x(t)$ is increasing (if $f\big(x(t),u(t)\big)>0$) or decreasing (if 
$f\big(x(t),u(t)\big)<0$). As long as $t\mapsto f(x(t),u(t))$ does not vanish on $(t_1,t_2)$,  
one can use the change of variable $s=x$ in the following integrals and obtain 
expressions of the elapsed time between instants $t_1$ and $t_2$
\begin{equation}\label{EQ_temps_reparam}
t_2-t_1=\int_{t_1}^{t_2}\!\!\!\mathrm{d}\tau=\int_{x(t_1)}^{x(t_2)}\!\!\!\frac{\mathrm{d}s}{f(s,u)},
\end{equation}
of the covered length between instants $t_1$ and $t_2$
\begin{equation}\label{EQ_moy_reparam}
\int_{t_1}^{t_2}\!\!\!x(\tau)\mathrm{d}\tau=\int_{x(t_1)}^{x(t_2)}
\!\!\!\frac{s\mathrm{d}s}{f(s,u)},
\end{equation}
and of the total energy consumption between instants $t_1$ and $t_2$:
\begin{equation}\label{EQ_cons_reparam}
\int_{t_1}^{t_2}\!\!\!h(x(\tau),u(\tau))\mathrm{d}\tau=
\int_{x(t_1)}^{x(t_2)}\!
\frac{h(s,u)\mathrm{d}s}{f(s,u)}.
\end{equation}
Indeed, the number of switches $N(u)$ involved in (\ref{EQ_def_cout}) is zero in formula (\ref{EQ_cons_reparam}) since $u$ is constant on $[t_1,t_2]$.
%\begin{equation}\label{EQ_cons_reparam}
%\int_{t_1}^{t_2}h(x(\tau),u)\mathrm{d}\tau=
%\left \{ \begin{array}{ll}
%\int_{x(t_1)}^{x(t_2)}
%\frac{h(s,1)\mathrm{d}s}{f(x_1,s,t,u)} & \mbox{ if }u=1\\
%0 & \mbox{ if }u=0 
%\end{array}\right.
%\end{equation}
From formulas (\ref{EQ_temps_reparam}), (\ref{EQ_moy_reparam}) and (\ref{EQ_cons_reparam}), one can notice that the speed average and the energy consumption are continuously differentiable with respect to the piece-wise constant control $u$, as long as the acceleration does not vanish. 
In particular, for any given constant $u$, for any given $V_0$, $V_1$ and any $w$ such that $x\mapsto f(x,u)$ does not vanish on the convex hull of 
$\{V_0,V_1,V_1+w \}$, the difference between the time needed to reach speed $V_1+w$ from $V_0$  and the 
time needed to reach speed $V_1$ from $V_0$ expresses, at first order, as
\begin{equation}\label{EQ_variation_temps}
\int_{V_0}^{V_1+w} \!\!\!\frac{\mathrm{d}s}{f(s,u)} -\int_{V_0}^{V_1} \!\!\!\frac{\mathrm{d}s}{f(s,u)}=\frac{w}{f(V_1,u)}+o_{w\to 0} (w).
\end{equation} 
Similarly, one gets the first order variation of the covered length
\begin{equation}\label{EQ_variation_distance}
\int_{V_0}^{V_1+w} \!\!\!\frac{s\mathrm{d}s}{f(s,u)} -\int_{V_0}^{V_1} \!\!\!\frac{s\mathrm{d}s}{f(s,u)}=\frac{wV_1}{f(V_1,u)}+o_{w\to 0} (w),
\end{equation} 
and of the energy consumption
\begin{equation}\label{EQ_variation_cout}
\int_{V_0}^{V_1+w} \!\!\frac{h(s,u)\mathrm{d}s}{f(s,u)} -\int_{V_0}^{V_1} \!\!\frac{h(s,u) 
\mathrm{d}s}{f(s,u)}\!=\!\frac{h(V_1,u)w}{f(V_1,u)}+o_{w\to 0} (w).
\end{equation} 
The above estimates are crucial for the variational arguments we will use in the proof of Propositions \ref{PRO_sol_autonome_non_periodique} and  \ref{PRO_structure_minimale}.

\subsection{Physical assumptions}
We list below some mathematical assumptions, that correspond to reasonable physical assumptions.
%Some mathematical assumptions, listed below, that correspond to reasonable physical assumptions, are necessary to ensure that Problem \ref{PROB_general_autonomous} has solutions. 
%Reasonable Physical assumptions make these solutions relevant. We now take advantage of the reparametrizations introduced in Section \ref{SEC_preliminaries_reparam} to give mathematically usable interpretations of the physical assumptions we need.
\begin{assumption}\label{ASS_physique_simplified}
\begin{enumerate}
 \item \label{ASS_acceleration_not_zero_simplified} The external actions like the slope of the track and the wind are not strong enough to prevent the car to go forward (maybe very slowly) when the engine is on: there exists  $V^\ast>0$ such $f(x,1)>0$ if $x < V^\ast$, $f(V^\ast,1)=0$ and $f(x,1)<0$ if $x>V^\ast$;
 \item \label{ASS_effective_engine}
  The (possibly negative) acceleration is larger when engine is switched on:  for every $x$, $f(x,0)<f(x,1)$;
\item \label{ASS_V_inf_simplified} When engine is off, the velocity tends to a (possibly negative) limit speed:
there exists $V_\ast$ such that $f(x,0)<0$ if $x>V_\ast$, $f(x,0)>0$ if $x<V_\ast$ and $f(V_\ast,0)=0$.  
\item Energy consumption is zero when the engine is off, positive when engine is on: for every $x$, $h(x,1)>h(x,0)=0$;\label{ASS_consumption_zero_motor_off}
\item Energy consumption increases with speed: $x\mapsto h(x,1)$ is non-decreasing; \label{ASS_cost_not_decreasing_simplified}
%\item (Generic transversality condition)for every $(x_1,t)$ in $[0,L]\times[0,T]$,  $\frac{d}{dx_2}f(x_1,V_\ast(x_,t),t,0)\neq 0$ and $\frac{d}{dx_2}f(x_1,V^\ast(x_,t),t,1)\neq 0$ ;
\item\label{ASS_switch_simplified}  The switching cost is not so large that the best strategy on long term is to keep full speed: all the integrals below converge and  $$\int_{V_\ast}^{V^\ast} \frac{h(V^\ast,1)-h(s,1)}{f(s,1)}\mathrm{d}s  +\alpha <~~~~~~~~ ~~~~~~~~$$
$$
~~~~~~~~\frac{h(V^\ast,1)}{V^\ast-V_\ast} \left ( 
\int_{V^\ast}^{V_\ast} \frac{s-V_\ast}{f(s,0)}\mathrm{d}s 
+ \int_{V_\ast}^{V^\ast}  \frac{s-V^\ast}{f(s,1)}\mathrm{d}s  \right );
$$
\item\label{ASS_conv_simplified} The energetic cost does not grow at the same rate than the engine efficiency when the velocity rises: 
the function  $F:x \mapsto \frac{h(x,1) f(x,0)}{f(x,1)-f(x,0)}$ is  either strictly concave or strictly convex. 
\end{enumerate}
\end{assumption} 

\begin{figure}
\begin{center}
\resizebox{10cm}{!}{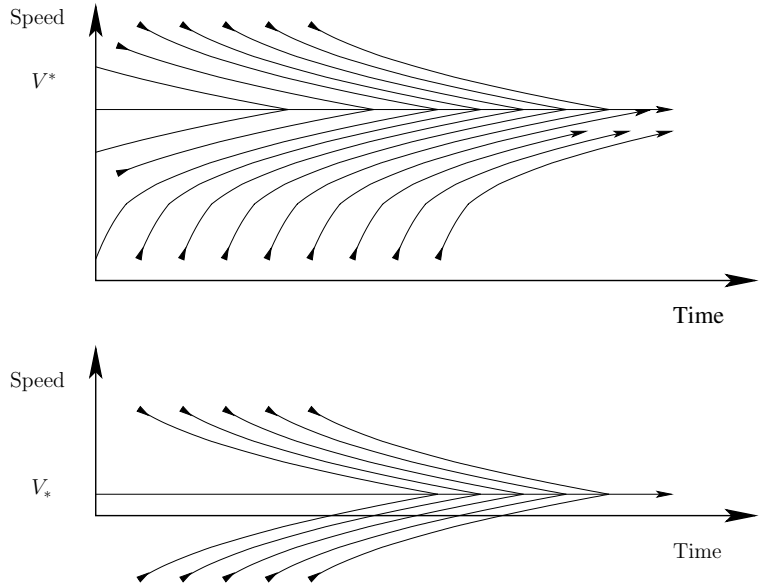}
\caption{Time plot of the dynamics with constant control $u=1$ (top) and $u=0$ (bottom). All the trajectories tend to $V^\ast$ (top) or $V_\ast$ (bottom) as the time tends to infinity. }\label{FIG_phase_portrait}  
\end{center}
\end{figure}

Assumptions \ref{ASS_physique_simplified}.\ref{ASS_acceleration_not_zero_simplified},   
 \ref{ASS_physique_simplified}.\ref{ASS_effective_engine}, and \ref{ASS_physique_simplified}.\ref{ASS_V_inf_simplified} are summarized in Figure~\ref{FIG_phase_portrait}.

The various terms involved in Assumption \ref{ASS_physique_simplified}.\ref{ASS_switch_simplified} can be interpreted in terms of velocity and accelerations. In particular, the first term 
%\begin{eqnarray*}
%\lefteqn{\int_{V_\ast}^{V^\ast} \frac{h(V^\ast,1)-h(s,1)}{f(s,1)}\mathrm{d}s =}\\
%&~~~~~~~~~~~~~~~~~~~& h(V^\ast,1)\int_{V_\ast}^{V^\ast}\frac{\mathrm{d}s}{f(s,1)} - \int_{V_\ast}^{V^\ast} \frac{h(s,1)}{f(s,1)}\mathrm{d}s
%\end{eqnarray*} 
$
\int_{V_\ast}^{V^\ast} \frac{h(V^\ast,1)-h(s,1)}{f(s,1)}\mathrm{d}s $ is equal to $ h(V^\ast,1)\int_{V_\ast}^{V^\ast}\frac{\mathrm{d}s}{f(s,1)} - \int_{V_\ast}^{V^\ast} \frac{h(s,1)}{f(s,1)}\mathrm{d}s,
$
i.e., the opposite of the difference between the energy needed to accelerate from $V_\ast$ to $V^\ast$ and the energy used by the vehicule running at full speed $V^\ast$ during the same amount of time $\int_{V_\ast}^{V^\ast}\frac{\mathrm{d}s}{f(s,1)}$.
%The other terms have similar interpretations that will be made clear in the proof of Lemma \ref{LEM_croissance_infini}.

\subsection{Oscillating structure of the solutions of Problem \ref{PROB_general_autonomous}}

\begin{proposition} \label{PRO_sol_autonome_non_periodique}
  Under Assumptions 
  \ref{ASS_maths_simplified} and \ref{ASS_physique_simplified}, 
either  Problem \ref{PROB_general_autonomous} has no solution or all of the solutions of Problem \ref{PROB_general_autonomous} have the 
following structure: with the possible exceptions of one initial and one final 
acceleration/deceleration phase, the speed of the vehicle oscillates 
between two values without intermediate local extrema.
\end{proposition}

 \begin{figure}
 \begin{center}
 \includegraphics[scale=0.6]{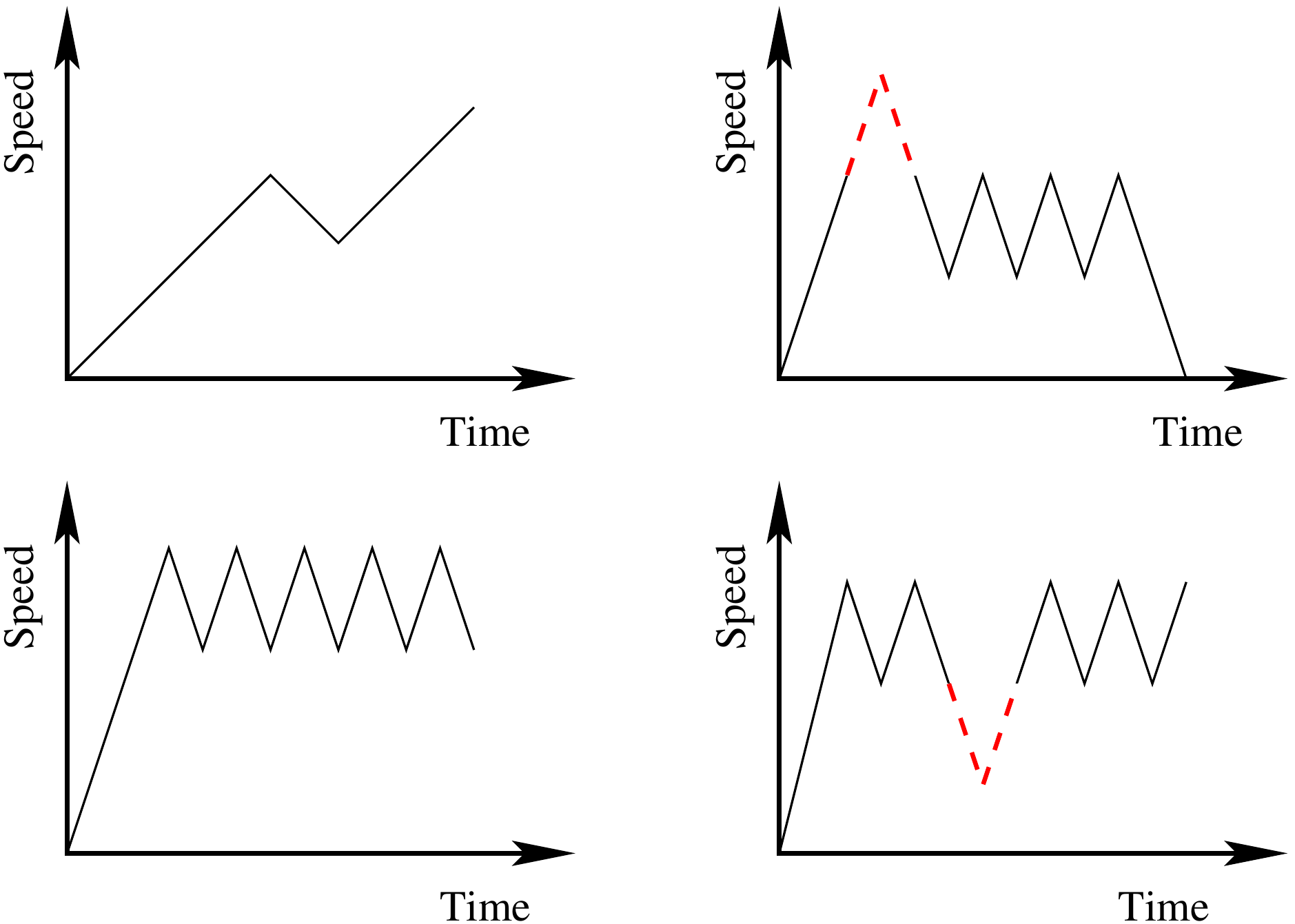}
 \caption{Representation of the speed with respect to the time. In this example, the initial speed $x_0$ is assumed to be zero. Both pictures in the left column represent possible optimal solutions to Problem \ref{PROB_general_autonomous}. Both pictures on the right cannot be possible optimal
 solutions of Problem \ref{PROB_general_autonomous}
  because of the red dashed  parts.}
 \label{FIG_optima_possibles}
\end{center}
\end{figure} 
 \begin{remark}
Proposition \ref{PRO_sol_autonome_non_periodique} does not state that Problem \ref{PROB_general_autonomous} has a solution. It only states that any solution of Problem \ref{PROB_general_autonomous} has an oscillating structure (see Figure \ref{FIG_optima_possibles}), which is no surprise. Indeed, if it is possible to write a Maximum Principle for this optimization problem, then the co-state and the associated switching function will be parametrized by the velocity $x$. As a general rule, there is a Maximum Principle for every optimization problem  but a precise statement may require space and technicalities (see in particular the very stimulating discussion of \cite[Section VI]{832814}). Notice also that the lack of backward uniqueness of trajectories prevents a direct use of the classical proofs of the Maximum Principle that can be found for instance in \cite{MR2062547}.
\end{remark}

The rest of this section is dedicated to the proof of Proposition \ref{PRO_sol_autonome_non_periodique}. The core of the proof is the following variational argument (Lemma \ref{LEM_argument_variation}) whose proof does not require the convexity guaranteed by Assumption \ref{ASS_physique_simplified}.\ref{ASS_conv_simplified}.
\begin{lemma}\label{LEM_argument_variation}
Under Assumptions 
  \ref{ASS_maths_simplified} and 
  \ref{ASS_physique_simplified}.\ref{ASS_acceleration_not_zero_simplified}
to   \ref{ASS_physique_simplified}.\ref{ASS_cost_not_decreasing_simplified}, 
 the solutions  (if any) of Problem \ref{PROB_general_autonomous} have the 
following structure: with the possible exceptions of one initial and one final 
acceleration/deceleration phase, the speed of the vehicle oscillates 
between local extrema $V_1,V_2,\ldots,V_p$ in $[V_\ast,V^\ast]$ and the matrix
$$
\begin{pmatrix}
1 & 1 & \ldots & 1 \\
V_1 & V_2 & \ldots & V_p\\
\frac{h(V_1) f(V_1,0)}{f(V_1,0)-f(V_1,1)} & \frac{h(V_2) f(V_2,0)}{f(V_2,0)-f(V_2,1)} & \ldots & \frac{h(V_p) f(V_p,0)}{f(V_p,0)-f(V_p,1)}
\end{pmatrix}
$$
has rank at most two.
\end{lemma}

\textbf{Proof of Lemma \ref{LEM_argument_variation}}:
Let $u:[0,T]\to \{0,1\}$ be a piecewise constant  optimal solution of Problem \ref{PROB_general_autonomous}. The restriction to $(0,T)$ of the associated trajectory $x:[0,T]\to \mathbf{R}$ admits successive local extrema $V_1,V_2,\ldots,V_p$ in $[V_\ast, V^\ast]$  at points $T_{V_1},T_{V_2},\ldots,T_{V_p}$ where $u$ has discontinuities, successive global maxima $V^\ast$ on  some intervals of length $T_1,T_2,\ldots,T_l$ and successive global minima $V_\ast$ on  some intervals of length $T'_1,T'_2,\ldots,T'_{l'}$. 
%All these extrema may be entangled, $x(T_{V_1})=V_1$, $x(T_{V_2})=V_2$ and $x$   being constant equal to $V^\ast$ on an interval of length $T_1$ included in $[T_{V_1},T_{V_2}]$. We do neither exclude the case when $x$ oscillates between $V_\ast$ and $V^\ast$ without intermediate local extrema, nor the case when $x$ never reaches $V_\ast$ or $V^\ast$ nor the case when $T_i=0$ or $T'_i=0$ for some integers $i$.  
We insist on the fact that $V_1,V_2,\ldots,V_p$ belong to the \emph{closed} interval $[V_\ast, V^\ast]$. That is, the same control $u$ may be associated with several different representations. For instance, a control $u$ that  switches from 1 to 0 exactly when the speed reaches $V^\ast$ can be represented either with some $T_{V_i}$ with $V_i=V^\ast$, or with some $T_i=0$ (meaning that the maximal speed is constant during a time equal to zero). Note also that we only consider extrema of $x$ on the \emph{open} interval $(0,T)$ (i.e., we do not consider the initial and final acceleration/deceleration if any). 
\begin{figure}
\begin{center}
\resizebox{15cm}{!}{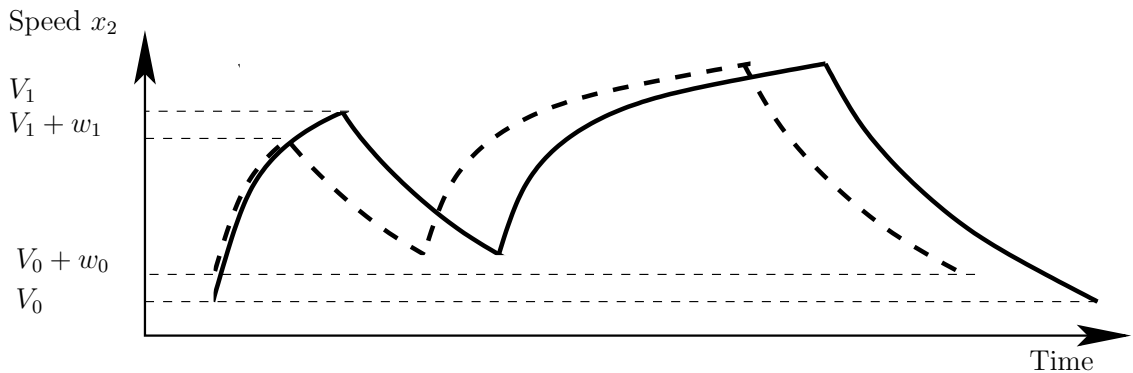}
\caption{First order variation argument. The reference trajectory (solid line), with extremal speeds $V_j$, $1\leq j\leq 4$, is compared with the close neighbor (dashed line) with extremal speeds $V_j+w_j$, $1\leq j\leq 4$. In this example, $w_2=w_3=0$.}  
\end{center}
\end{figure}
 We proceed with a classical argument of calculus 
of variation and introduce $u^{w,\tau}$,
 an other admissible control for average speed $V$, 
 close to $u$ in norm $L^1$,
which admits local minima $V_j+w_j$, local maxima $V_j+w_j$ and (if applicable for the reference trajectory) global extrema $V_\ast$ and $V^\ast$ on  some intervals of length $T_1+\tau_1,T_2+\tau_2,\ldots,T_l+\tau_l$. If 
$T_i=0$ or $V_i=V_\ast$, we impose $\tau_i\geq 0$. If $V_i=V^\ast$, we impose $w_i \leq 0$.  

From (\ref{EQ_variation_temps}), the invariance of total time (race duration equals $T$) implies
\begin{equation}\label{EQ_contrainte_temps}
\sum_{j=1}^{p}  w_j \left (\frac{1}{f(V_j,1)}-\frac{1}{f(V_j,0)} \right )+ \sum_{j=1}^l \tau_j+ \sum_{j=1}^{l'} \tau'_j =0.
\end{equation}
From (\ref{EQ_variation_distance}), the invariance of the average speed (covered length is equal to $L$) implies
\begin{equation}\label{EQ_contrainte_vitesse}
\sum_{j=1}^{p}   
 \left (\frac{w_j V_j}{f(V_j,1)}-\frac{w_j V_j}{f(V_j,0)} \right )
+V^\ast \sum_{j=1}^l \tau_j + V_\ast \sum_{j=1}^{l'} \tau'_j =0.
\end{equation}
 The optimality of $u$  implies the invariance at first order of the cost. Hence, from (\ref{EQ_variation_cout}), 
 $$
\sum_{j=0}^{p}  \frac{w_jh(V_j,1)}{f(V_j,1)} +h(V^\ast,1) \sum_{j=1}^l \tau_j  + h(V_\ast,0) \sum_{j=1}^{l'} \tau'_j =0.
 $$
 From Assumption \ref{ASS_physique_simplified}.\ref{ASS_consumption_zero_motor_off}, 
 $$
\sum_{j=0}^{p}  \frac{w_jh(V_j,1)}{f(V_j,1)} +h(V^\ast,1) \sum_{j=1}^l \tau_j   =0,
 $$ 
 or, defining $\tilde{w}_j=w_j \left (\frac{1}{f(V_j,1)}-\frac{1}{f(V_j,0)} \right )$, $\Theta=\sum_{j=1}^l \tau_j $ and $\Theta'=\sum_{j=1}^{l'} \tau'_j$ 
\begin{equation}\label{EQ_contrainte_consommation}
\sum_{j=0}^{p}  \tilde{w}_j \frac{h(V_j,1)f(V_j,0)}{f(V_j,0)-f(V_j,1)} +h(V^\ast\!\!,1) \Theta  =0.
\end{equation}

Considering the left hand sides of  equations (\ref{EQ_contrainte_temps}), (\ref{EQ_contrainte_vitesse}) and (\ref{EQ_contrainte_consommation}) as linear forms in $(\tilde{w}_j,\Theta,\Theta')$, we infer that the matrix
$$
\begin{pmatrix}
1 & \ldots & 1 & 1 & 1\\
V_1 &  \ldots & V_p & V^\ast & V_\ast\\
\frac{h(V_1,1) f(V_1,0)}{f(V_1,0)-f(V_1,1)} & \ldots & \frac{h(V_p,1) f(V_p,0)}{f(V_p,0)-f(V_p,1)} & h(V^\ast,1) & 0
\end{pmatrix}
$$
has rank at most two (the last two columns only appears when $l\neq 0$ or $l'\neq 0$).  Finally, notice that 
$$
 h(V^\ast,1)=  \frac{h(V^\ast,1) f(V^\ast,0)}{f(V^\ast,0)-f(V^\ast,1)}
$$
and 
$$
 h(V_\ast,0)=  \frac{h(V_\ast,1) f(V_\ast,0)}{f(V_\ast,0)-f(V_\ast,1)}=0
$$
to complete the proof of Lemma \ref{LEM_argument_variation}. \hfill ~$\blacksquare$

\textbf{Proof of Proposition \ref{PRO_sol_autonome_non_periodique}:}
We proceed by contradiction and assume that the speed admits $p$ different local extrema. If $p=2$,  there is nothing to prove. Else,
we consider the first three columns of the matrix obtained in Lemma \ref{LEM_argument_variation}. Up to permutation of the indices, we may assume without loss of generality that $V_1<V_3<V_2$. Find $\mu,\nu >0$, $\mu+\nu=1$ such that $V_3=\mu V_1+\nu V_2$. Since the matrix
$$
\begin{pmatrix}
1 & 1  & 1=\mu.1 +\nu.1 \\
V_1 & V_2  & V_3=\mu V_1 +\nu V_2\\
\frac{h(V_1,1) f(V_1,0)}{f(V_1,0)-f(V_1,1)} & \frac{h(V_2,1) f(V_2,0)}{f(V_2,0)-f(V_2,1)}  & \frac{h(V_3,1) f(V_3,0)}{f(V_3,0)-f(V_3,1)}
\end{pmatrix}
$$
has rank at most 2, its determinant is zero. Hence,
\begin{eqnarray*}
\lefteqn{\frac{h(V_3,1) f(V_3,0)}{f(V_3,0)-f(V_3,1)}=}\\
&~~~~& \mu \frac{h(V_1,1) f(V_1,0)}{f(V_1,0)-f(V_1,1)} + \nu  \frac{h(V_2,1) f(V_2,0)}{f(V_2,0)-f(V_2,1)}.
\end{eqnarray*}
Hence, the mapping $F:x\mapsto \frac{h(x,1) f(x,0)}{f(x,1)-f(x,0)}$ is not strictly convex or concave, which is in contradiction  with 
Assumption \ref{ASS_physique_simplified}.\ref{ASS_conv_simplified}. 
\hfill ~$\blacksquare$

\section{Autonomous dynamics with infinite time horizon}\label{SEC_simplified_periodic}

\subsection{Statement of the optimal problem}

In this subsection, we consider Problem \ref{PROB_general_autonomous} with a control restricted to a periodic control. Precisely, for any given $V>0$, we define 
the \emph{admissible trajectories with average speed $V$} as the curves 
$t\in[0,T_1] \mapsto (x(t),u(t)) \in \mathbf{R}\times \{0,1\}$ for some suitable $T_1>0$, with $x$ absolutely 
continuous and $u$ piecewise constant, satisfying (\ref{EQ_dynamics_simplified}) and 
such that $x(0)=x(T_1)=V$, $u(0)=u(T_1)$ 
and $\frac{1}{T_1}\int_0^{T_1}x(t)\mathrm{d}t=V$. 
%Notice that we consider only trajectories for which the 
%final velocity $x(T_1)$ is equal to the initial velocity $x(0)$ and the final value of 
%the control $u(T_1)$ is equal to initial one $u(0)$.
In particular, the control $u$ has an 
even number of discontinuities in $[0,T_1]$.

A piecewise constant function $u:[0,T_1]\to \{0,1\}$ is an \emph{admissible control for 
average speed $V$} if there exists a function $x:[0,T_1]\to \mathbf{R}$ such that $(x,u)$
is an admissible trajectory with average speed $V$.

For $T_1\geq 0$ and $k$ in $\mathbf{N}$, we denote with $\mathcal{U}^k_V(T_1)$  
the set of admissible controls for 
average speed $V$ defined on $[0,T_1]$ with at most $k$ discontinuities, 
and we define the sets
$\mathcal{U}_V(T_1)=\cup_{k \in \mathbf{N}}\mathcal{U}^k_V(T_1)$, 
$\mathcal{U}^k_V=\cup_{T_1>0}\mathcal{U}^k_V(T_1)$ and
 $\mathcal{U}_V=\cup_{T_1>0}\mathcal{U}_V(T_1)=\cup_{k \in \mathbf{N}}\mathcal{U}^k_V$.  

For any subset $I$ of $\mathbf{R}$, we denote respectively with 
$\mathcal{U}_V^{k,I}(T_1)$, $\mathcal{U}_V^{I}(T_1)$ and $\mathcal{U}_V^{I}$ the functions of $\mathcal{U}_V^{k}(T_1)$, $\mathcal{U}_V(T_1)$ and $\mathcal{U}_V$ where $x$ takes value in $I$.

Finally, we associate, with every admissible control for average speed $V$ in 
$\mathcal{U}_V(T_1)$, the cost
$$C^a(u)=\frac{1}{T_1}C(u)$$
which can be seen as 
the average energy consumption to run the distance $V T_1$ in time $T_1$. Let us consider the following problem.

\begin{problem}\label{PROB_autonomous_periodic} Let $V>0$ be given. Find the minimum, if any, of $C^a$ on  $\mathcal{U}_V^I$.
\end{problem}

Problem \ref{PROB_autonomous_periodic} may be seen as an infinite time horizon version of Problem \ref{PROB_general_autonomous}. Indeed, no final time is given. We only aim at finding the lowest possible energy consumption for a long distance run at average speed $V$.

\subsection{Main Result}
The main result, proved below, stipulates that the minimum of $C^a$ is reached by controls with two 
discontinuities. In other words, the most energy efficient driving strategy is to let the 
speed $x$ periodically oscillate between two values $V^a$ and $V^b$, with no 
intermediate local extremum. Precisely, one has the following result (recall that $V^\ast$ and $V_\ast$ are defined in Assumption \ref{ASS_physique_simplified}.\ref{ASS_acceleration_not_zero_simplified}): 
\begin{proposition}\label{PRO_structure_minimale}
Under Assumptions \ref{ASS_maths_simplified} and \ref{ASS_physique_simplified},  for every closed interval $I$  contained in 
 $[V_\ast,V^\ast]$ and every  interior point $V$ of $I$, $C^a$ admits a minimum on $\mathcal{U}_V^I$, and 
$\displaystyle{\inf_{\mathcal{U}_V^I}C^a\!\!=\min_{\mathcal{U}^{2,I}_V}C^a}$. \end{proposition}

 The result of Proposition \ref{PRO_structure_minimale} is twofold: first, there exists a minimum of $C^a$ and second, this minimum is reached with periodic controls with two discontinuities for a period.  The second part is a direct application of Proposition \ref{PRO_sol_autonome_non_periodique}. The only point to check is the existence of solutions, which is done in Section \ref{SEC_proof}. 

 \subsection{Proof}
 \label{SEC_proof}
 The proof requires several intermediate results and is given in the following section.
 \subsubsection{Existence of solutions to Problem \ref{PROB_autonomous_periodic} for a given number of switches} 
 \label{SEC_proof_exist}
%To begin with, we assume that there is no constraint on the state that is  $I=[V_\ast,V^\ast]$. 
\begin{lemma}\label{LEM_existence_min}
For every even $k$ in $\mathbf{N}$ and $T_1>0$, 
$C^a$ (and hence $C$) admits a 
minimum $m_V^k(T_1)$ on $\mathcal{U}^k_V(T_1)$.
\end{lemma}
\textbf{Proof:}
 This is a consequence of the compactness of 
 $\mathcal{U}^k_V(T_1)$ endowed with the $L^1$ norm and of the continuity of $C$ when the number of switches is constant equal to~$k$.\hfill  ~ $\blacksquare$
 %%%%%%%%%%%%%%%%%%%
 \subsubsection{Cost estimates for control with small periods}
   %%%%%%%%%%%%%%%%
\begin{lemma}\label{LEM_explosion_0}
 For every  $k$ in $2\mathbf{N}$,
$ \lim_{T_1 \to 0}m_V^k(T_1) =+\infty$
\end{lemma} 
\textbf{Proof:} We first prove by contradiction that  $\mathcal{U}^k_V(T_1)$ does not contain any continuous 
function.  Indeed, assume that $u$ in $\mathcal{U}^k_V(T_1)$ is continuous. Then,  the 
associated trajectory  $x$ is decreasing (if $u=0$) or increasing (if $u= 1$). 
Since $u$ belongs to $\mathcal{U}^k_V(T_1)$, $x(0)=x(T_1)$. Hence $x$ is constant. 
By Assumptions \ref{ASS_physique_simplified}.\ref{ASS_acceleration_not_zero_simplified} and \ref{ASS_physique_simplified}.\ref{ASS_V_inf_simplified}, 
$x= V^\ast$ or $x= V_\ast$. 
Hence, the average of $x$ is not $V$ (which belongs to $(V_\ast,V^\ast)$).
 This is a contradiction with the fact that $u$ belongs to $\mathcal{U}^k_V(T_1)$. Hence $u$ has at least two discontinuities and 
$N(u)\geq 1$. Thus, $C^a(u)>\alpha/T_1$. \hfill ~~$\blacksquare$

\begin{remark}
As $\alpha$ tends to zero, the interval between two consecutive switching times tends to 
zero. Hence, the velocity tends to a constant $V$, which is the optimal trajectory of the 
relaxed problem when $u$ can take value in $[0,1]$ with $f(x_1,x_2,u)=u f(x_1,x_2,1)+(1-
u)f(x_1,x_2,0)$. 
\end{remark}

\subsubsection{Oscillating structures of optimal trajectories with a given period}
\begin{lemma}\label{LEM_structure_T_fini}
For every $k$ in $2\mathbf{N}$, for every positive $T_a$, there exists $p$ in 
$\mathbf{N}$ such that $m^k_V(T_a)=m^2_V(T_a/p)$ 
\end{lemma}
\textbf{Proof:}
This is nothing but a restatement of Proposition~\ref{PRO_sol_autonome_non_periodique} in the periodic case.
\hfill ~ $\blacksquare$
\subsubsection{Cost estimates for controls with large periods}

 \begin{lemma}\label{LEM_croissance_infini}
For every $T_1$ large enough, 
$m_V^2(T_1)\leq \liminf_{T_2\to +\infty} m_V^2(T_2)$.
\end{lemma}
\textbf{Proof:} We do the proof by computing an expansion of $m_V^2(T_2)$
as $T_2$ tends to infinity. Let $T_2$ be given, and $u$ a minimizing control in  $m_V^2(T_2)$. We distinguish among several cases, depending on the regularity of $f$.

\begin{figure}
\begin{center}
\resizebox{8cm}{8cm}{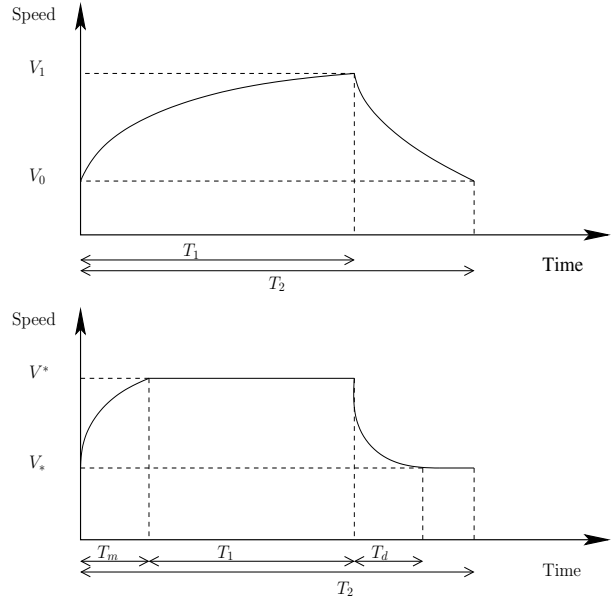}
\caption{Large period optimal controls. For a control in $\mathcal{U}^2_V(T_2)$, we define $T_1$ as the acceleration duration in the case where trajectories are backward unique near $V_1$ (top) and as the time where $x=V^\ast$ when backward uniqueness does not hold (bottom). The ratio $T_1/T_2$ is the same at first order in both cases.} 
\label{FIG-Tm}
\end{center}
\end{figure}

In the case of local backward uniqueness of the trajectories of (\ref{EQ_dynamics_simplified}) in  neighborhoods of $V^\ast$ and $V_\ast$ (see Fig. \ref{FIG-Tm} top), we may assume, up to a time shift,  that $u(t)=0$ for $t\in(0,T_1)$ and $u(t)=0$ for $t\in (T_1,T_2)$. The associated trajectory reaches its minimum $V_0$ at $t=0$ and its maximum $V_1$ at $T_1$. 
From (\ref{EQ_temps_reparam}), we get that
$$
\int_{V_0}^{V_1} \frac{s}{f(s,1)}\mathrm{d}s=\int_{V_0}^{V_1} \frac{s-V^\ast}{f(s,1)}\mathrm{d}s + V^\ast T_1 \mbox{ and }
$$
$$
\int_{V_1}^{V_0} \frac{s}{f(s,0)}\mathrm{d}s=\int_{V_1}^{V_0} \frac{s-V_\ast}{f(s,0)}\mathrm{d}s + V_\ast (T_2-T_1).
$$
From \eqref{EQ_moy_reparam}, we get that
\begin{eqnarray*}
T_2 V&=&  \int_{V_0}^{V_1} \frac{s}{f(s,1)}\mathrm{d}s + \int_{V_1}^{V_0} \frac{s}{f(s,0)}\mathrm{d}s
\end{eqnarray*}
Hence, 
\begin{eqnarray}
\lefteqn{V^\ast \frac{T_1}{T_2}+V_\ast \left( 1-\frac{T_1}{T_2} \right )
+\frac{1}{T_2} \int_{V_0}^{V_1}\frac{s-V^\ast}{f(s,1)}\mathrm{d}s}\nonumber\\
&~~~~~~~~~~~~~~~~~~~~~~~~~~~~~~~~~~~& +
\frac{1}{T_2} \int_{V_1}^{V_0}\frac{s-V_\ast}{f(s,0)}\mathrm{d}s=V
\end{eqnarray}
that implies 
$\lim_{T_2\to \infty}T_1=+\infty$. Hence, $\lim_{T_2\to \infty}V_0=V_\ast$ and 
$\lim_{T_2\to \infty}V_1=V^\ast$. Thus, there exists  a function $\epsilon$ such that $\lim_{+\infty}\epsilon=0$ and
\begin{eqnarray}
\lefteqn{\frac{T_1}{T_2}=-\frac{1}{T_2 (V^\ast-V_\ast)} \left ( \int_{V_\ast}^{V^\ast}\!\!\!\! \frac{s-V^\ast}{f(s,1)}\mathrm{d}s + \int_{V^\ast}^{V_\ast}\!\!\! \frac{s-V_\ast}{f(s,0)}\mathrm{d}s \right )}\nonumber \\
&~~~~~~~~~~~~~~~~~~~~~~~~~~~~~~~~~~~~&+\frac{V-V_\ast}{V^\ast-V_\ast}  +\frac{\epsilon(T_2)}{T_2}.\label{EQ_est_T1}
\end{eqnarray}
We infer from (\ref{EQ_temps_reparam}) and (\ref{EQ_est_T1}) that
\begin{eqnarray}
\lefteqn{C^a= \frac{\alpha}{T_2}+\frac{1}{T_2}\int_{V_0}^ {V_1}\frac{h(s,1)}{f(s,1)}\mathrm{d}s}\nonumber \\
&=&  \frac{\alpha}{T_2}+\frac{1}{T_2}\left( \int_{V_0}^{V_1} \frac{h(s,1)-h(V^\ast,1)}{f(s,1)}\mathrm{d}s +h(V^\ast,1) T_1 \right)\nonumber \\
&=& h(V^\ast,1) \frac{V-V_\ast}{V^\ast-V_\ast}+\frac{1}{T_2}\left \lbrack \alpha +\int_{V_\ast}^ {V^\ast}\!\!\!\frac{h(s,1)-h(V^\ast,1)}{f(s,1)}\mathrm{d}s  \right .  \nonumber \\
&&\left . -\left ( \int_{V_\ast}^ {V^\ast}\frac{s-V^\ast}{f(s,1)}\mathrm{d}s + \int_{V^\ast}^ {V_\ast}\frac{s-V_\ast}{f(s,0)}\mathrm{d}s \right ) \frac{h(V^\ast,1)}{V^\ast-V_\ast} \right \rbrack  \nonumber\\
&&
+\frac{\epsilon(T_2)}{T_2}
\label{EQ_est_Ca}
\end{eqnarray}

In the case where the trajectories of (\ref{EQ_dynamics_simplified}) reach $V^\ast$ and $V_\ast$ in finite time, we denote (see Fig. \ref{FIG-Tm} bottom)
$
T_m=\int_{V_\ast}^{V^\ast}\frac{\mathrm{d}s}{f(s,1)}$ { and }
$T_d=\int_{V_\ast}^{V^\ast}\frac{\mathrm{d}s}{f(s,0)}.
$ 
Hence, 
\begin{equation}\label{EQ_est_T1_sature}
\frac{T_1}{T_2}(V^\ast-V_\ast)=V-V_\ast+(Tm+T_d) \frac{V}{T_2} 
\end{equation}
Plugging (\ref{EQ_est_T1_sature}) in the definition of $C^a$, we find that the estimate (\ref{EQ_est_Ca}) is still valid (with $\epsilon=0$) when the trajectories of (\ref{EQ_dynamics_simplified}) are not backward uniquely defined near $V^\ast$ and $V_\ast$. 
The mixed case (the trajectories are locally backward unique  near $V^\ast$ and reach $V_\ast$ in finite time, or vice versa) can be treated as the two cases above and is omitted. 

In conclusion, the estimate (\ref{EQ_est_Ca}) is valid in any case and Lemma \ref{LEM_croissance_infini}  
follows from Assumption \ref{ASS_physique_simplified}.\ref{ASS_switch_simplified}.
\hfill ~$\blacksquare$

%The second point is to notice that $C(u)\geq \alpha \frac{TV(u)}{T_1}$. Hence, for every 
%$T_1$, $\inf_k m_V^k(T_1)=\min_k m_V^k(T_1)$. 

From Lemmas \ref{LEM_explosion_0} and \ref{LEM_croissance_infini}, we deduce 
%(see Figure  \ref{FIG_allure_m_V}) 
that the continuous 
 mapping   $T \mapsto m_V^k(T)$ admits a minimum on $(0,+\infty)$. We denote this quantity with $m_V^k:=\min_{T_1>0} m_V^k(T_1)$.

\subsubsection{Completion of the proof of Proposition \ref{PRO_structure_minimale}}
 By definition,
$\inf_{\mathcal{U}_V}C^a = \inf_{T>0}\inf_{k \in 2\mathbf{N}} m^k_V(T)$.
By  Lemma \ref{LEM_existence_min}, each optimization problem for a given period $T_1$ admits a solution, and by Lemma \ref{LEM_structure_T_fini}, 
the associated control can be chosen with two discontinuities
 ${\inf_{\mathcal{U}_V}C^a = \inf_{T>0}  m^2_V(T),}$
 this last quantity being equal to $m^2_V$ by Lemmas \ref{LEM_explosion_0} and \ref{LEM_croissance_infini}. In particular, the infimum  $\inf_{\mathcal{U}_V}C^a $ is reached and is indeed a minimum. ~~\hfill~ $\blacksquare$

\subsection{Practical determination of the upper and lower speed limits} \label{SEC_algorithm}

Based on the knowledge of the average target speed $V$, the dynamics $f$ and the cost function $h$, 
the search for the upper and lower speeds of the optimal solution whose existence 
is asserted by Proposition \ref{PRO_structure_minimale} amounts to solve two 1D 
continuous optimization problems. A variety of methods are available in the 
literature,  especially when some regularity assumptions are done on $f$ and $h$. 
We present below a naive yet efficient grid method.

\textbf{Step 1: candidate selection} Select $p$ candidates $V^a_1<V^a_2<\ldots <V^a_p<V$ for the lower speed limit.\\
\textbf{Step 2: search for upper limit} For every $i=1\ldots p$,  find (by dichotomy), the upper limit $V^b_i$ ensuring that the average speed is $V$ (in the case where $V^b_i=V^\ast$, compute also the time where the velocity is constant).\\
%, i.e., in the case $(V^a_i,V^b_i) \subset (V_\ast,V^\ast)$, solve for $V^b_i$:
%$$
%\frac{\int_{V^a_i}^{V^b_i} s\left(\frac{1}{f(s,1)}-\frac{1}{f(s,0)} \right)\mathrm{d}s}{\int_{V^a_i}^{V^b_i} \left(\frac{1}{f(s,1)}-\frac{1}{f(s,0)} \right)\mathrm{d}s}
%= V. 
%$$
\textbf{Step 3: cost for each candidate} For every $i=1\ldots p$, compute the cost $C^a$ associated with a trajectory oscillating from $V^a_i$ to $V^b_i$ (or staying at $V^\ast$ for the time computed in Step 2).\\
%, equal to
%$$
%\frac{\alpha + \int_{V^a_i}^{V^b_i} \frac{sh(s)}{f(s,1)}\mathrm{d}s}%{\int_{V^a_i}^{V^b_i} \left(\frac{1}{f(s,1)}-\frac{1}{f(s,0)} \right)\mathrm{d}s}.
%$$
\textbf{Step 4: conclusion} Pick the pair $(V^a_i,V^b_i)$ with the lowest cost.

In our implementation (see Section \ref{SEC_implementation}), we mainly used $p=4$,  $V^a_i=V-\frac{5-i}{2}$ m.s$^{-1}$, $i=1..4$. 

\section{Non autonomous problem in finite time}\label{SEC_Framework}
\subsection{Modeling}\label{SEC_Modeling}

For a real vehicle, the autonomous dynamics described in Section \ref{SEC_Autonome_simplified} is too restrictive. Various types of forces act on the vehicle. Some of them, like the solid 
friction forces applied on the wheels or axes,  are roughly constant.
Some others, such as the gravitation, depends on the slope of track. 
%Some forces depend only on the velocity (such as ). 
Some of the forces also depend on many (partly unknown) external factors, including the actual weather.
This is the case for
%the additional torque provided by the solar cells depending among others on the clouds opacity or for
the aerodynamic drag depending on   the 
velocity, the orientation of the track with respect to the wind and the strength of the  wind.

If we consider both the position and the velocity of the vehicle denoted respectively $x_1$ and $x_2$, the dynamics turns into:
\begin{equation}\label{EQ_dynamic}
\left \{ 
\begin{array}{lcl}
\dot{x}_1&=&x_2\\
\dot{x}_2&=&f(x_1,x_2,t,u)
\end{array} \right.
\end{equation}
with $f:[0,L]\times [0,+\infty) \times [0,T] \times \{0,1\} \to \mathbf{R}$ 
a  continuous function. The control $u:[0,T]\to \{0,1\}$ is a piecewise constant function 
accounting for the state of the engine (off when $u=0$, on when $u=1$).

The steering possibilities of the car being limited at high speed, we introduce a (smooth)  
function $x_1\mapsto V^s(x_1)$ representing the maximal safety speed at point $x_1$ and we will 
only consider solutions of system (\ref{EQ_dynamic}) satisfying the state constraint $0\leq x_2(t)
\leq V^s(x_1(t))$ for every $t$.

As in Section \ref{SEC_problem_statement_simplified}, with every admissible control $u$ and absolutely continuous solution $(x_1,x_2)$ of (\ref{EQ_dynamic}), we associate the energy cost
\begin{equation}\label{EQ_def_cout_non_autom}
C(u)=\int_0^T h(x_2(t),u(t))\mathrm{d}t+ \alpha\,N(u)
\end{equation}
where $h:[0,+\infty)\times \{0,1\} \to [0,+\infty)$ is a $C^1$ function, $\alpha$ is a 
positive real constant and $N$ denotes the number of discontinuities of $u$.

\begin{problem} We aim to find a piecewise constant control $u$ that
steers $(x_1,x_2)$, subject to the velocity constraints $x_2(t)\in V^s(x_1(t))$ for every $t$ and the dynamic (\ref{EQ_dynamic}),  from 
$(0,0)$ at time $t=0$ 
to $(L,x_2(T))$ in time less than $T$ and minimizes cost $C$ \eqref{EQ_def_cout_non_autom}.
\end{problem}

We do  the following assumptions.
\begin{assumption}\label{ASS_maths}
\begin{enumerate}
\item \label{ASS_pcw_continuous} The track and weather conditions are piecewise continuous: there exist two subdivisions  $0=x_1^0<x_1^1<\ldots<x_1^{p_1}=L$ and
$0=t^0<t^1<\ldots <t^{p_t}=T$ of $[0,L]$ and $[0,T]$ such that for every $u$ in $\{0,1\}$, $(x_1,x_2,t)\mapsto f(x_1,x_2,t,u)$ is continuous
on $[x_1^j,x_1^{j+1})\times \mathbf{R}\times [t^l,t^{l+1})$ for every $j<p_1$ and $l<p_t$;
\item \label{ASS_forward_uniqueness} For $u$ constant in $\{0,1\}$, for every $j<p_1$ and $l<p_t$, for every  $(x_1(t_0),x_2(t_0),t_0)$ in $[x_1^j,x_1^{j+1})\times \mathbf{R}\times [t^l,t^{l+1})$,
 the Cauchy problem (\ref{EQ_dynamic}) with initial condition $(x_1(t_0),x_2(t_0),t_0)$ at time $t_0$ admits a locally unique solution in positive time;
\item \label{ASS_local_simplifiable} For every time $t^\ast$ in $[0,T]$, for every position on the track $x_1^\ast$ in $[0,L]$, the functions $(x,u)\mapsto f(x_1^\ast, x, t^\ast, u)$ and $(v,u)\mapsto h(v,u)$ satisfy Assumptions \ref{ASS_maths_simplified} and \ref{ASS_physique_simplified}. 
\end{enumerate}
\end{assumption} 
Assumption \ref{ASS_maths}.\ref{ASS_local_simplifiable} means that, around $x_1^\ast$ in $[0,L]$, we can locally  replace the actual vehicle evolving on the actual irregular track with varying weather conditions by a fictional vehicle evolving on an infinite fictional track with constant slope or wind velocity  (at $x_1^\ast$).

%\ref{ASS_ass}.\ref{ASS_conv}   
%The physical meaning of Assumptions \ref{ASS_ass}.\ref{ASS_switch} and 
%\ref{ASS_ass}.\ref{ASS_conv} will become clear in Section \ref{SEC_proof}.

Let $u:[0,T_1]\to \{0,1\}$ be piecewise constant. A pair 
$(x_1,x_2):[0,T_1]\to \mathbf{R}^2$ is a  \emph{solution} of 
(\ref{EQ_dynamic}) if $(x_1,x_2)$ is  continuous, piecewise differentiable and 
satisfies (\ref{EQ_dynamic}) for  almost every time $t$ in $[0,T_1]$.

As before, the existence (but not the uniqueness) of such solutions to Cauchy problem (\ref{EQ_dynamic}) with initial condition $(x_1(t_0),x_2(t_0))$ at time $t_0$ is a consequence of Assumption \ref{ASS_maths}.\ref{ASS_pcw_continuous} and Peano theorem (see for instance 
\cite{MR0171038}, Theorem 2.1 page 10).

%\section{Adaptative strategy and analysis}\label{SEC_Analysis}

\subsection{Description of the strategy}\label{SEC_strategy}
The strategy we propose can be split in two parts: first, we  identify the current dynamics  and compute the current target and then we use the results of Section~\ref{SEC_Autonome_simplified} during a small time interval (typically three seconds in our case), after which we proceed to a new identification, and a new optimization.
  
Precisely, we split the time interval $[0,T]$ into small intervals $[0,T]=\bigcup_{n=0}^N [n t_a, (n+1)t_a]$. 

Assume that current time is equal to $t=n t_a$ for some $n$, and the vehicle is at position $x_1(n t_a)$ with velocity $x_2(n t_a)$. Identify the current dynamics \eqref{EQ_dynamic}
%(which is an approximation of (\ref{EQ_dynamics_simplified}))
 with $x_1^\ast=x_1(n t_a)$ and $t^\ast=n t_a$ and
compute the target average speed $V=\frac{L-x_2(n t_a)}
{T-n t_a}$. If $V< V^s(x_1(t_a))$,  find the minimum $m_V^2$ and the corresponding trajectory $(x_2,u^\ast)$  of (\ref{EQ_dynamics_simplified}), with the constraint $x_2\leq V^s(x_1(t_a))$. This can be done, for instance, with Step~1 to Step~4 detailed in Section~\ref{SEC_algorithm}.  The optimal speed $x_2$ oscillates between two values $V^a(t_a)$ and $V^b(t_a)$. 

If $V^b(t_a)> V^s(x_1(t_a)$, keep $V^b(t_a)=V^s(x_1(t_a))$ and $V^a(t_a)=V^s(x_1(t_a))-\delta$  
where $\delta$ is a small enough constant (for our application, $\delta =0.5$ m.s$^-1$).

The control strategy $\mathcal{S}$ is to repeat the following until $t=(n+1)t_a$:
\begin{enumerate}
\item Keep $u=1$ until $x_2(t)$ reaches $V^b(n t_a)$. 
\item As soon as $x_2(t)$ reaches $V^b(n t_a)$, set $u=0$.
\item Keep $u=0$ until $x_2(t)=V^a(t)$, then switch to $u=1$.
\item Come back to (1), unless $x_1(t) \geq L$.
\end{enumerate}

\subsection{Well-posedness}
The only point to check to guarantee the well-posedness of the control scheme $\mathcal{S}$ 
is that the interval between two consecutive engine switches 
cannot tend to zero. Indeed, let $u_t$ be an admissible optimal control for average speed $V$ 
computed at time $t$. We denote the period of $u_t$ with $T_t$ and the associated cost with $C^a_t$ (index $t$ is the time at which the control has been computed). By Lemma \ref{LEM_explosion_0}, $C^a(u_t)>\alpha /T_1$, hence $T_t > \alpha / C^a(u_t)$. The cost $t\mapsto C^a_t$ is continuous and the interval $[0,T]$ is compact, hence $\inf_{0\leq t \leq T} C^a_t>0$.   
  And thus, $\inf_{t\in [0,T]} T_t >0$, and 
 the switching points cannot accumulate (no Zeno phenomena).

\subsection{Constraints fulfillment and robustness}\label{SEC_stability}

The proposed strategy happens to be extremely robust in practice. This robustness may be explained by two reasons of different natures. 

The first reason for the robustness of the strategy is the continuous actualization of the target average speed.

The second reason for robustness is the fact that we use a control on the average of the velocity. 
%Assume that we made an error in the identification of the dynamics, and we have overestimated (rep. underestimated) the slope of the track. Then, the vehicle will need less (resp. more) time than expected to reach the upper speed limit starting form  the lower limit. But the average speed during the process will be almost what was expected. The aim of this Section \ref{SEC_stability} is to  formalize this intuitive result.
We first state   two abstract results (Lemma  \ref{LEM_robustesse_variation_premier_ordre} and Proposition \ref{PRO_invariance_moyenne_premier_ordre} below).

\begin{lemma}\label{LEM_robustesse_variation_premier_ordre}
Let $V_a<V_b$ be given, denote with $\mathcal C$ the set of not vanishing continuous function on $[V_a,V_b]$. For $g$ in $\mathcal C$, denote 
$\displaystyle{
T(g)=\int_{V_a}^{V_b} \frac{\mathrm{d}s}{g(s)}}$, $\displaystyle{L(g)=\int_{V_a}^{V_b} \frac{s\mathrm{d}s}{g(s)}}$ and the mapping $\mathcal F$ defined on $\mathcal C$ by 
$
\mathcal{F}:g \mapsto L(g)/T(g) \in \mathbf{R}. 
$ .
Then, for every $g$ in $\mathcal C$ and every continuous $\delta g$ such that $\|\frac{\delta g}{g} \|_{L^\infty}<1$, 
\begin{eqnarray}
\lefteqn{\mathcal{F}(g+\delta g)-\mathcal{F}(g)}\nonumber \\
&= &\frac{1}{T(g+\delta g)} \sum_{n=1}^{\infty} (-1)^n \! \int_{V_a}^{V_b} \!\!\frac{s-\mathcal{F}(g)}{g(s)} \left ( \frac{\delta g(s)}{g(s)} \right )^n \!\!\!\mathrm{d}s.\label{EQ_variation_robustesse_premier_ordre}
\end{eqnarray}
\end{lemma}

\textbf{Proof:}
Standard expansions show that 
$$\displaystyle{L(g+\delta g)=\sum_{n\geq 0} (-1)^n \int_{V_a}^{V_b} \frac{s (\delta g)^n}{g^{n+1}(s)}\mathrm{d}s},$$ and
$\displaystyle{T(g+\delta g)=\sum_{n\geq 0} (-1)^n\int_{V_a}^{V_b} \frac{ (\delta g)^n}{g^{n+1}(s)}\mathrm{d}s}$, for $\delta g$ small enough.
Conclusion follows with basic calculus. \hfill~ $\blacksquare$

\begin{proposition} \label{PRO_invariance_moyenne_premier_ordre}
If $\delta g/g$ is constant, then $\mathcal{F}(g+\delta g)=\mathcal{F}(g)$.
\end{proposition}
\textbf{Proof:} Apply (\ref{EQ_variation_robustesse_premier_ordre}) and notice that 
$\displaystyle{
\int_{V_a}^{V_b} \frac{s-\mathcal{F}(g)}{g(s)}\mathrm{d}s =}$
$L(g) - \mathcal{F}(g) T(g)=0$.
\hfill~ $\blacksquare$

% In practice, we apply Lemma \ref{LEM_robustesse_variation_premier_ordre} and Remark \ref{REM_invariance_moyenne_premier_ordre} with $g:s \mapsto f(s,1)$.

Let us comment the statements of 
Lemma \ref{LEM_robustesse_variation_premier_ordre} and 
Proposition \ref{PRO_invariance_moyenne_premier_ordre}. 
If $f$ satisfies Assumptions \ref{ASS_maths_simplified} 
and \ref{ASS_physique_simplified} 
and one choses $g:s\mapsto f(s,1)$, the number $\mathcal{F}(g)$ is the 
average speeds when one speed up from $V_a$ to $V_b$ with the dynamics $f$. If 
the dynamics $s\mapsto f(s,1)$ has been incorrectly identified and so is 
$\tilde{g}:s\mapsto f(s,1)+\delta f(s)$, 
Lemma \ref{LEM_robustesse_variation_premier_ordre} gives an estimate of the 
difference between  the expected average speed $\mathcal{F}(\tilde{g})$ and 
the actual average speed  $\mathcal{F}({g})$. The point of 
Proposition \ref{PRO_invariance_moyenne_premier_ordre} is that, if 
$(\delta f(s))/f(s,1)$ is constant, then the change of average speed is zero. Hence the difference $\mathcal{F}(g+\delta g)-\mathcal{F}(g)$ is due to the \emph{variance} of $\frac{\delta g}{g}$ which is usually small with respect to $\frac{\delta g}{g}$. The same computation is obviously valid for decelerations 
if one considers $g:s\mapsto f(s,0)$, and generically in the non-autonomous case if one
 considers $g:s\mapsto f(x_1(s),s,t(s),u)$ where $x_1(s)$ and $t(s)$ are respectively the position on the track and the time where the vehicle reaches speed $s$.

\section{Effective implementation}\label{SEC_implementation}

The strategy proposed in Section \ref{SEC_strategy} has been implemented in an official competition in 2014.

\subsection{Vir'Volt 2 prototype}

The prototype (including the pilot) has a total weight of about 90 kg, is driven by a 200 W 
electric DC motor powered by a 23 V battery. 
The control scheme presented in this note is implemented on a microcontroller (ref dsPIC33ep512mu810 from Microchip\textsuperscript{TM}), with a 140 MHz internal clock. 
Standard 32 kHz oscillators, bike velocity sensors and
GPS receiver provide (respectively) the current time, the velocity and the position needed for the 
algorithm. 
%The total cost of the power management unit including the board, the wires and the sensors is about 1000 euros.

A simplified dynamics for Vir'Volt 2 is given by
\begin{eqnarray}
\lefteqn{f(x_1,x_2,t,u)=%\nonumber \\&&
-a(x_2\!-v(x_1,t))^2-\!c\,\mathrm{sign}(x_2)-\!g\sin\theta_{x_1}}\nonumber \\
&&+\left \{\begin{array}{llll}
f_1 &\mbox{ if }u=1 \\ %&\mbox{ and } &x_2<v_{lim}\\
0&\mbox{ if }u=0 %&\mbox{ or } &x_2 > v_{lim}
\end{array} \right. \label{EQ_dynamique_Virevolt}
\end{eqnarray}
where $v(x_1,t)$ is the speed of the wind at point $x_1$ at time $t$, $\theta_{x_1}$ is the angle of the track with the horizontal plane at point $x_1$, 
%$v_{lim}=15.8$ m.s$^{-1}$ is the maximal speed for which power transmission to the wheels is possible (given by the engine nominal speed $\eta=5680$ rpm, the transmission chain with ratio 9 and wheel diameter $d_{wheel}=0.48 $m),
 $f_1$ is the traction force by mass unit when the engine is on,  $g$ is the gravitational acceleration, $a$ accounts for the aerodynamics drag, and $c$ accounts for the solid friction in the vehicle.
% $f_1=0.20$m.s$^{-2}$ is the traction force by mass unit when the engine is on,  $g=9.81$  m.s$^{-2}$ is the gravitational acceleration, $a=6.10^{-4}$ m$^{-1}$ accounts for the aerodynamics drag, and $c=3. 10^{-2}$ m.s$^{-2}$ accounts for the solid friction in the vehicle.
\begin{figure}
\centering
\begin{tabular}{cllc}
Constant & \multicolumn{2}{c}{Value} & Comment\\
\hline
$a$ & $6.10^{-4}$ &m$^{-1}$ & identified\\
$c$ & $3.10^{-2}$ &m.s$^{-2}$ & identified\\
$g$ & $9.81$  &m.s$^{-2}$ & tabulated \\
$f_1$& $0.20$ &m.s$^{-2}$ & tabulated \\
$m$& 93 &kg & measured\\
\hline
\end{tabular}
\caption{Numerical value of the constants appearing in the dynamics (\ref{EQ_dynamique_Virevolt}).}
\end{figure}

%\begin{figure}[hbtp]
%\centering
%\includegraphics[angle=90,scale=0.4]{Dessins/Deceleration.pdf}
%\caption{An example of deceleration record  (a practically flat portion of the track during a free practice in Rotterdam 2014) of the prototype used for identification of the dynamics.The cross marks are obtained from an optic sensor measuring the position of the front wheel of the vehicle. The solid blue line is the best fit obtained letting $c$ change in the dynamics (\ref{EQ_dynamique_Virevolt}). This long time (18 seconds) deceleration is mainly used for off-line identification of $c$. On-line identification of the track slope $\theta$ is done with much shorter (3 seconds) decelerations.}  
%\label{FIG_Essais_Rotterdam} 
%\end{figure}

The traction force $f_1$ is tabulated from the engine data sheet.
%The constant $a$ has been identified after multiple (a few hundreds) decelerations with engine turned off on an indoor ($v=0$) flat horizontal track ($\theta=0$). The constant $c$ depends on the track and has been identified from data taken from free practice session in Rotterdam (65 decelerations, see Figure \ref{FIG_Essais_Rotterdam}).

There are several ways to consider the energetic cost. The power transmitted to the wheels is the product
$
h_1(x_2,1)=x_2\cdot m \cdot f_1$. %of the speed, the mass and the acceleration.  It grows linearly with the speed. 
The energy taken from the battery is always larger than $h_1$ and depends on the design of the power conversion block. In the following, we consider the case where the power consumption is roughly independent from speed,  $h_2(x_2,1)=23 \cdot 7=161$ W for $x_2> 5 \mbox{ m.s}^{-1}$.
Notice that the dynamics $f$ completed with $h_1$ or $h_2$ satisfies in both cases Assumption \ref{ASS_maths} at any time.

For a flat horizontal track without wind, with $\alpha=10$ J and a target average speed of 7 m/s, the optimal strategy consists in letting the speed oscillate  between 6.1  m/s and 7.94 m/s. The corresponding energetic cost for a track of 16.5 km (similar to the editions 2013 and 2014 of the European Shell Eco-Marathon) is 104 189 J. This computation is given for reference only since it neglects the fact that the vehicle start from velocity zero (kinetic energy of the vehicle at 7 m/s: about 2.5 kJ), the influence of the wind, and the elevation of the actual track which is not perfectly horizontal. 

\subsection{Effective results}

Computation of the control law (including identification of the current dynamics from the 
sensors data) is possible with an average power of 10 mW.  

The official results of the EMT prototype measured in the 2014 edition of the Shell Eco-Marathon are presented in Figure \ref{table_results} with the first result of the 2013 
edition.  The 
2013 prototype is basically the same as the 2014 one without embedded control (the driver 
directly switches the engine on or off). These results are directly comparable since they 
all were obtained without solar cells (EMT record on the Rotterdam track is 90211 J and was obtained with manual velocity command during a particularly sunny test in May 2013).

%\begin{figure*}[!t]
%\centering
%\subfloat[Case I]{\includegraphics[width=2.5in]{VitessesRun1.pdf}%
%\hfil
%\subfloat[Case II]{\includegraphics[width=2.5in]{VitessesRun1.pdf}%
%\label{fig_second_case}}
%\caption{Simulation results.}
%\label{fig_sim}
%\end{figure*}

% An example of a floating table. Note that, for IEEE style tables, the 
% \caption command should come BEFORE the table. Table text will default to
% \footnotesize as IEEE normally uses this smaller font for tables.
% The \label must come after \caption as always.
%
\begin{figure}[!t]
\renewcommand{\arraystretch}{1.2}
 %if using array.sty, it might be a good idea to tweak the value of
 %\extrarowheight as needed to properly center the text within the cells
%\centering
% Some packages, such as MDW tools, offer better commands for making tables
% than the plain LaTeX2e tabular which is used here.
\begin{center}
\begin{tabular}{c c c c}
\hline
& 2013 & \multicolumn{2}{c}{2014}\\
%\cline{3-4}
& 1st run & 1st run & 2nd run \\
\hline
Run duration \hfill ~ & & &\\
~~~~~expected   & 35m 00s & 36m 00s &38m 30s\\
~~~~~realized & 36m 01s & 36m 10s &38m 31s\\
%\hline
Consumption & 112 742 J & 120 150 J & 107 946 J\\
\hline
\end{tabular}
\caption{Comparisons of the official runs results for the EMT prototype 
at Shell Eco-Marathon
2013 (manual command) and 2014 (automatic command and strategy of Section \ref{SEC_strategy}).}
\label{table_results}
\end{center}
\end{figure}

\section{Conclusions}

A real-time adaptive driving command strategy for low consumption vehicle has been 
presented. Real life implementation on a competition prototype showed good 
agreement between the theoretical analysis and the effective results. The consumption 
performances are comparable with performances obtained by highly trained human drivers. All 
the computations are done on board, in real time, with a low cost micro controller and limited energy consumption.    

Among the various advantages of the automatic command, the pilots underlined the dramatic 
safety improvement (no need to concentrate on the velocity anymore) and the robustness of 
the command that allows effective recovering of optimal timing after traffic perturbations.\\

From a more theoretical point of view, it has been provided an explicit hybrid control strategy in a form of a piecewise constant control (on/off) for low regular dynamics, in particular not Lipschitz continuous. Furthermore, a switching cost has been considered. Because of the heavy computational burden, the use of a classical PMP is redhibitory. This contribution acts an an efficient alternative from this perspective. Optimality of the control in the case of autonomous systems and the robustness in the case of non autonomous systems (with time-varying disturbances) are shown.  
   
%Our assumptions of non-predictability of the characteristics of the track are
% certainly too conservative. The slope may vary in space but is obviously 
% constant in time and can 
% be incorporated in the 
% optimization process.  Modeling of the weather influence is more subtle yet is not out of reach. Indeed, the circuit is roughly circular with a small diameter 
% (about 250 m), hence the strength and direction of the wind can basically be recovered 
% at any point from one single measure. One should keep in mind that the resulting optimization process has to be done with low cost embedded devices. A balance has to be found between a more precise dynamics and the limited on board computational resources.

\section*{Ackowledgements}                               % Place acknowledgements
This work has been financially supported by CPER MISN and the F\'ed\'eration Charles Hermite.  
It is a pleasure for the authors to thank  the ``Eco Motion Team'' in charge of the Vir'Volt prototype development, especially G\'erard D\'echenaud, Jean-Claude Sivault, Mathieu Rupp  and Sandra Rossi.% here.

\bibliographystyle{plain}
\bibliography{References}

\end{document}